# Functional principal components analysis via penalized rank one approximation

**Jianhua Z. Huang**[*]

*Department of Statistics*
*Texas A&M University*
*College Station, TX 77843, USA*
*e-mail:* jianhua@stat.tamu.edu

**Haipeng Shen**[†]

*Department of Statistics and Operations Research*
*University of North Carolina*
*Chapel Hill, NC 27599, USA*
*e-mail:* haipeng@email.unc.edu

and

**Andreas Buja**

*Department of Statistics*
*The Wharton School*
*University of Pennsylvania*
*Philadelphia, PA 19104, USA*
*e-mail:* buja.at.wharton@gmail.com

**Abstract:** Two existing approaches to functional principal components analysis (FPCA) are due to Rice and Silverman (1991) and Silverman (1996), both based on maximizing variance but introducing penalization in different ways. In this article we propose an alternative approach to FPCA using penalized rank one approximation to the data matrix. Our contributions are four-fold: (1) by considering invariance under scale transformation of the measurements, the new formulation sheds light on how regularization should be performed for FPCA and suggests an efficient power algorithm for computation; (2) it naturally incorporates spline smoothing of discretized functional data; (3) the connection with smoothing splines also facilitates construction of cross-validation or generalized cross-validation criteria for smoothing parameter selection that allows efficient computation; (4) different smoothing parameters are permitted for different FPCs. The methodology is illustrated with a real data example and a simulation.

---

[*]Supported in part by by NSF grant DMS-0606580 and NCI grant CA57030.
[†]Supported in part by NSF grant DMS-0606577.

## Contents

## 1. Introduction

Principal components analysis is a key technique for unsupervised functional data analysis Ramsay and Silverman (2005) where the goal is to decompose variation in a two-way data table. Two different approaches to smoothed functional principal components analysis (FPCA) were proposed by Rice and Silverman (1991) and Silverman (1996), both of which we briefly describe in what follows. Let $X(\cdot)$ denote a random function which we assume for now can be observed repeatedly and as a whole, without discretization; in addition, let $\alpha$ be a smoothing penalty parameter. To find the $j$-th principal component weight function $\gamma_j(\cdot)$, the Rice-Silverman approach maximizes

$$\frac{\operatorname{var}(\int \gamma\,X) - \alpha \int \gamma''^2}{\int \gamma^2} \tag{1}$$

subject to the constraint $\int \gamma \hat\gamma_k = 0$ for $k < j$, where $\hat\gamma_k$ is the estimated $k$-th principal component function, while the Silverman approach maximizes

$$\frac{\operatorname{var}(\int \gamma\,X)}{\int \gamma^2 + \alpha \int \gamma''^2}. \tag{2}$$

subject to $\int \gamma \hat\gamma_k + \alpha \int \gamma'' \hat\gamma_k'' = 0$ for $k < j$. Both approaches impose smoothness on principal components using roughness penalty ideas, but they differ in the way they incorporate the penalty. The variance $\operatorname{var}(\int \gamma\,X)$ in (1) and (2) needs to be estimated from a random sample of realizations of $X(\cdot)$.

Following Hotelling (1933), maximizing variance of a standardized linear combination of variables is the standard textbook treatment of principal components analysis (PCA, e.g., Mardia et al., 1979). A different approach is by way of fitting low rank approximations to the data matrix (Pearson, 1901), and this approach is intimately connected to the singular value decomposition (Eckart and Young, 1936). See Jolliffe (2002) for further discussion and references. In this article we apply the roughness penalty idea to a low rank approach to FPCA. As we will see, there exist difficulties in introducing roughness penalties, and a guiding principle in navigating these difficulties will be certain invariance under scale transformations. Scale invariance considerations will not only allow us to select the proper form of penalty, but also to re-derive and thereby justify the regularized variance criterion of Silverman (1996). However, penalized low rank approximation combined with invariance principles amounts to more than a novel justification for an existing approach; it has several methodological advantages over existing regularized variance approaches.

First, our approach yields a power algorithm that is an efficient variant of the power algorithm (e.g., Appendix A of Jolliffe, 2002) for calculating eigenvectors. Second, spline smoothing of discretized data is naturally built into our method which therefore gains the theoretical advantages of smoothing splines: our method applies directly to discretized data, and the estimated FPCs are solutions of an optimization problem defined on a function space. More importantly, the connection of our method to spline smoothing helps us develop computationally efficient cross-validation (CV) and generalized cross-validation (GCV) criteria for selecting smoothing parameters. This development fills a gap in the FPCA literature as is shown by a comparison: both Rice and Silverman (1991) and Silverman (1996) assume the whole curve is available through either interpolation or smoothing, while (Ramsay and Silverman, 2005, Chapters 8 & 9) represents functional data using a basis expansion prior to performing PCA. Finally, our method allows different principal components to have different smoothing parameters, which is a beneficial flexibility that is not shared by the method of Silverman (1996).

In this paper we focus on functional data that are sampled on a common grid across subjects — a typical setting of functional data analysis. Sometimes, the functions may be irregularly or even sparsely sampled, as often occurs in biomedical longitudinal studies. FPCA methods for sparsely sampled functional data or longitudinal data have been developed, among others, by James (2000), Rice (2001), and Yao et al. (2005a,b).

The rest of the article is organized as follows. Section 2 presents the new method for extracting principal components for finite dimensional vector spaces; Section 5 extends the method to function space. Section 3 describes the power algorithm. Section 4 discusses the criteria for smoothing parameter selection and also provides the derivations. Section 6 gives some numerical results. Finally, the Appendix contains some technical details.

## 2. A new framework on functional principal components

Since functional data are usually observed discretely, this section presents our method in terms of discretized data. Unlike standard multivariate analysis, we will take into account the intrinsic functional structure of the data with general Ridge penalties, which in the functional case will be smoothness penalties. Recovery of the whole curve of principal component weight functions involves spline interpolation and will be discussed in Section 5. In what follows we focus on the first principal component; subsequent principal components can be extracted sequentially by removing preceding components.

### 2.1. Minimizing reconstruction error

Consider a collection or sample of functional data, observed or recorded at a common set of discrete observation points $t_1, \ldots, t_m$. Denote the underlying functions for the sample as $x_i(\cdot)$, $i = 1, \ldots, n$. The observed data are $x_{ij} = x_i(t_j)$, $i = 1, \ldots, n$, $j = 1, \ldots, m$. Let the $n \times m$ data matrix be $\mathbf{X} = (x_{ij})$. For simplicity and without loss of generality, suppose that $\mathbf{X}$ is column centered, that is, the sample mean of each column of $\mathbf{X}$ is zero.

Consider the problem of finding the best rank one approximation of $\mathbf{X}$. Any rank one matrix of size $n \times m$ can be written as $\mathbf{u}\mathbf{v}^T$, where $\mathbf{u} = (u_1, \ldots, u_n)^T$ and $\mathbf{v} = (v_1, \ldots, v_m)^T$. Let $\|\mathbf{X}\|_F$ denote the Frobenius norm of $\mathbf{X}$. The problem can be formally stated as minimizing with respect to $\mathbf{u}$ and $\mathbf{v}$ the following reconstruction error,

$$\|\mathbf{X} - \mathbf{u}\mathbf{v}^T\|_F^2 \triangleq \operatorname{tr}\{(\mathbf{X} - \mathbf{u}\mathbf{v}^T)(\mathbf{X} - \mathbf{u}\mathbf{v}^T)^T\} = \sum_{i=1}^n \sum_{j=1}^m (x_{ij} - u_i v_j)^2. \quad (3)$$

For a fixed $\mathbf{v}$, the $\mathbf{u}$ that minimizes (3) is $\mathbf{u} = \mathbf{X}\mathbf{v}/\mathbf{v}^T\mathbf{v}$. Plugging this $\mathbf{u}$ into (3) we find that the minimizing $\hat{\mathbf{v}}$ of (3) maximizes $\mathbf{v}^T\mathbf{X}^T\mathbf{X}\mathbf{v}/\mathbf{v}^T\mathbf{v}$. Since $\mathbf{X}^T\mathbf{X}$ is proportional to the sample variance-covariance matrix, $\mathbf{v}^T\mathbf{X}^T\mathbf{X}\mathbf{v}/\mathbf{v}^T\mathbf{v}$ is proportional to the sample variance of the data projected onto the direction of $\mathbf{v}$. Thus $\hat{\mathbf{v}}$ is, up to a scale factor, the first principal component loading (weight) vector for $\mathbf{X}$ according to standard multivariate analysis textbooks (e.g., Mardia et al., 1979). We also note that the minimizing $\hat{\mathbf{u}}$ and $\hat{\mathbf{v}}$ of (3) are the first singular vector pair of $\mathbf{X}$, again up to a scale factor.

The appeal of minimizing (3) is that it resembles a prediction problem. If $\mathbf{u}$ were observable, then (3) would be the least squares criterion for a multivariate regression. The difference is that the predictor variable $\mathbf{u}$ needs to be estimated here. Based on the above observation, we next discuss how to modify (3) to take into account the functional nature of data by smoothing with a roughness penalty. The connection with a prediction problem will be critical once again for deriving cross-validation criteria in Section 4.

### 2.2. Penalization and transformation invariance

Following the basic philosophy of functional data analysis, we assume that there is an underlying smooth function $\gamma(\cdot)$ such that $v_j = \gamma(t_j)$. If the domain of values $t$ is the real line, we can assume $t_i$ sorted. Because of the smoothness of $\gamma(\cdot)$, a pair of adjacent values $v_j$ and $v_{j+1}$ are necessarily tied to each other by correlation. Using the roughness penalty idea, it is natural to consider minimizing

$$\|\mathbf{X} - \mathbf{u}\mathbf{v}^T\|_F^2 + \alpha\,\mathbf{v}^T\mathbf{\Omega}\mathbf{v}, \qquad (4)$$

where $\mathbf{\Omega}$ is a non-negative definite roughness penalty matrix, and $\alpha > 0$ is a penalty parameter. The penalty matrix $\mathbf{\Omega}$ is chosen such that a larger value of $\mathbf{v}^T\mathbf{\Omega}\mathbf{v}$ is associated with a greater penalty on differences between adjacent values. For example, an intuitive choice of $\mathbf{\Omega}$ for equispaced $t_j$ based on second differences of $\mathbf{v}$ would be $\mathbf{v}^T\mathbf{\Omega}\mathbf{v} = \sum_{j=2}^{m-1}(v_{j+1} - 2v_j + v_{j-1})^2$. Construction of general penalty matrices that are suitable for both smoothing and interpolation will be discussed in Section 5 (see Theorem 1).

An immediate problem with (4) is that it is not invariant under the pair of scale transformations $\mathbf{u} \to c\mathbf{u}$ and $\mathbf{v} \to \mathbf{v}/c$. While the goodness-of-fit term $\|\mathbf{X} - \mathbf{u}\mathbf{v}^T\|_F^2$ does not change, the roughness penalty term $\mathbf{v}^T\mathbf{\Omega}\mathbf{v}$ changes. A simple possible fix is to make the penalty term unitless and to minimize

$$\|\mathbf{X} - \mathbf{u}\mathbf{v}^T\|_F^2 + \alpha\,\frac{\mathbf{v}^T\mathbf{\Omega}\mathbf{v}}{\mathbf{v}^T\mathbf{v}}. \qquad (5)$$

Fixing $\mathbf{v}$, the minimizing $\mathbf{u}$ is $\mathbf{u} = \mathbf{X}\mathbf{v}/\mathbf{v}^T\mathbf{v}$, which can be plugged back into (5) to show that minimizing (5) is equivalent to maximizing

$$\frac{\mathbf{v}^T\mathbf{X}^T\mathbf{X}\mathbf{v} - \alpha\,\mathbf{v}^T\mathbf{\Omega}\mathbf{v}}{\mathbf{v}^T\mathbf{v}}. \qquad (6)$$

This maximizing criterion is essentially the same as the Rice-Silverman criterion expressed in (1).

However, the criterion (5) has a defect: The optimization problem is not invariant under scale transformation of the measurements. Under the transformation $\mathbf{X} \to c\mathbf{X}$ and the corresponding transformation $\mathbf{u} \to c\mathbf{u}$, the goodness-of-fit term becomes $c^2\|\mathbf{X} - \mathbf{u}\mathbf{v}^T\|_F^2$, while the penalty term remains unchanged. The reason is that the goodness-of-fit term is not unitless but the penalty term is. The criterion (4) has the same defect because of a mismatch of units of the two terms involved. — This motivates us to consider minimizing

$$\|\mathbf{X} - \mathbf{u}\mathbf{v}^T\|_F^2 + \alpha\,\mathbf{u}^T\mathbf{u}\,\mathbf{v}^T\mathbf{\Omega}\mathbf{v}, \qquad (7)$$

where the two terms now have the same units. Optimization of (7) is invariant under scale transformation in the following sense: If $\hat{\mathbf{u}}$ and $\hat{\mathbf{v}}$ form its minimizer for the data matrix $\mathbf{X}$, then $\hat{\mathbf{u}}^* = c\hat{\mathbf{u}}$ and $\hat{\mathbf{v}}^* = \hat{\mathbf{v}}$ form its minimizer for the rescaled data matrix $\mathbf{X}^* = c\mathbf{X}$. In other words, the minimizing $\mathbf{v}$ of this criterion is the same before and after the scale transformation $\mathbf{X} \to c\mathbf{X}$. — Fixing $\mathbf{v}$, the minimizing $\mathbf{u}$ is $\mathbf{u} = \mathbf{X}\mathbf{v}/\mathbf{v}^T(\mathbf{I}_m + \alpha\,\mathbf{\Omega})\mathbf{v}$, which can be plugged back into (7) to show that minimizing (7) is equivalent to maximizing

$$\frac{\mathbf{v}^T\mathbf{X}^T\mathbf{X}\mathbf{v}}{\mathbf{v}^T\mathbf{v} + \alpha\,\mathbf{v}^T\mathbf{\Omega}\mathbf{v}}. \qquad (8)$$

This maximizing criterion is essentially the same as the Silverman criterion expressed in (2).

Use of the penalized reconstruction error combined with invariance considerations indicates that the Silverman proposal is preferable to the Rice-Silverman proposal for functional principal components analysis. Such insight would not be forthcoming by studying penalized variance alone. In addition, as we shall show below, a direct methodological advantage of using the penalized reconstruction error (7) is that it facilitates extension of CV- and GCV- type smoothing parameter selection criteria for spline smoothing to FPCA.

From now on we shall focus on the criterion (7) for extracting smoothed principal components. An extension of (7) for extracting the whole principal component weight function is given in Section 5. Since we extract the principal component functions sequentially, multiple smoothing parameters are allowed. This is different from Silverman (1996) where a single smoothing parameter is used for extracting all principal component functions. The benefit of the flexibility provided by using multiple smoothing parameter will be illustrated in Section 6 with simulated data.

**Remark 1.** Invariance of scale transformation of measurements is motivated by the consideration that results should remain the same under change of metric of the data. Such consideration of invariance might be meaningless if the data matrix would be standardized prior to analysis. However, there is no obvious way of standardization for functional data. If each column of the data matrix is standardized by dividing the corresponding sample standard deviation, for example, then the functional nature of the data will be lost. Standardizing the whole data matrix by the overall sample standard derivation is not a sound operation either, because the sample variance may vary from column to column.

## 3. The power algorithm

The criterion (7) can be minimized by alternating minimization of $\mathbf{u}$ and $\mathbf{v}$ in an iterative algorithm: Fixing $\mathbf{v}$, $\mathbf{u} = \mathbf{X}\mathbf{v}/\mathbf{v}^T(\mathbf{I}_m + \alpha\,\mathbf{\Omega})\mathbf{v}$; fixing $\mathbf{u}$, $\mathbf{v} = (\mathbf{I} + \alpha\,\mathbf{\Omega})^{-1}\mathbf{X}^T\mathbf{u}/\mathbf{u}^T\mathbf{u}$. Separating the scale constants, the algorithm is as follows:

1. Initialize $\mathbf{v}$.
2. Repeat until convergence:

    (a) $\mathbf{u} \leftarrow \mathbf{X}\mathbf{v}$,

    (b) $\mathbf{v} \leftarrow (\mathbf{I} + \alpha\,\mathbf{\Omega})^{-1}\mathbf{X}^T\mathbf{u}$,

    (c) $\mathbf{v} \leftarrow \mathbf{v}/\|\mathbf{v}\|$.

The initial $\mathbf{v}$ can be chosen to be the first right singular vector of $\mathbf{X}$. Step 2. (c) forces $\mathbf{v}$ to have norm 1, which is somewhat arbitrary and only meant for identifiability purposes between $\mathbf{u}$ and $\mathbf{v}$; any other normalization with different scale trade-offs between $\mathbf{u}$ and $\mathbf{v}$ would work, too. For example, an alternative restriction on $\mathbf{v}$ is $\mathbf{v}^T(\mathbf{I} + \alpha\,\mathbf{\Omega})\mathbf{v} = 1$, but its dependence on the smoothing parameter $\alpha$ renders it less convenient. When starting from the right singular vector of $\mathbf{X}$, this algorithm converges rather quickly, usually in a few iterations. The algorithm is a simple variant of the power algorithm for computing eigenvectors. The iteration in terms of $\mathbf{v}$ alone is $\mathbf{v} \leftarrow (\mathbf{I} + \alpha\,\mathbf{\Omega})^{-1}\mathbf{X}^T\mathbf{X}\mathbf{v}; \ \mathbf{v} \leftarrow \mathbf{v}/\|\mathbf{v}\|$. If there is no penalty ($\alpha = 0$), the algorithm is essentially the power algorithm for conventional PCA (e.g., page 409 of Jolliffe, 2002). As discussed below in Section 4, this iterative algorithm naturally facilitates the derivation of an explicit cross-validation criterion for smoothing parameter selection.

## 4. Choosing the smoothing parameter

### *4.1. Cross-validation and generalized cross-validation*

Our formulation of FPCA suggests a new method for selecting the smoothing parameter. Step 2. (b) of the iterative algorithm essentially smooths $\mathbf{X}^T\mathbf{u}$ with $\mathbf{S}(\alpha) = (\mathbf{I} + \alpha\,\mathbf{\Omega})^{-1}$ to update $\mathbf{v}$. This interpretation suggests that CV and GCV criteria for selecting the spline tuning parameter $\alpha$ (Green and Silverman, 1994) can be adopted to FPCA. The CV score is defined as

$$\mathrm{CV}(\alpha) = \frac{1}{m}\sum_{j=1}^{m} \frac{\left[\{(\mathbf{I}-\mathbf{S}(\alpha))(\mathbf{X}^T\mathbf{u})\}_{jj}\right]^2}{\left(1 - \{\mathbf{S}(\alpha)\}_{jj}\right)^2}, \qquad (9)$$

and the GCV score is defined as

$$\mathrm{GCV}(\alpha) = \frac{1}{m}\, \frac{\left\|(\mathbf{I}-\mathbf{S}(\alpha))(\mathbf{X}^T\mathbf{u})\right\|^2}{\left(1 - \frac{1}{m}\mathrm{tr}\{\mathbf{S}(\alpha)\}\right)^2}. \qquad (10)$$

In the implementation of either method, we can simply nest CV- or GCV-selection of $\alpha$ inside the loop, i.e., in Step 2. (b) of the algorithm.

The above motivation of CV and GCV is by analogy only. To give a formal justification, we show below in Section 4.2 that the CV score given above can indeed be derived from the basic idea of cross-validation, deletion of one data item at a time. Similarly, the GCV score can also be derived from the basic idea of GCV (Craven and Wahba, 1979). The technical difference of our derivations compared to those for smoothing splines is that we delete *one column* of $\mathbf{X}$ at a time, rather than *one matrix entry*.

The connection of column deletion in FPCA with point deletion in spline smoothing is clearly seen by considering an extreme case of (7). If $\mathbf{X}$ has only one row, denoted by $\mathbf{y}^T$, then $\mathbf{u} = u$ is a scalar. Requiring $\mathbf{u}$ to have norm 1 and fixing its sign for identifiability, it is necessary that $u = 1$. Then (7) becomes

$$\|\mathbf{y} - u\mathbf{v}\|^2 + \alpha\,u^2\mathbf{v}^T\mathbf{\Omega}\mathbf{v} = \|\mathbf{y} - \mathbf{v}\|^2 + \alpha\,\mathbf{v}^T\mathbf{\Omega}\mathbf{v},$$

which is exactly the penalized least squares criterion for smoothing splines. This connection motivates our column deletion procedure that we now elaborate.

### 4.2. Derivation of CV and GCV criteria

Minimization of (7) with $\mathbf{u}$ fixed can be considered as a ridge regression. Conditional on $\mathbf{u}$ with $\mathbf{v}$ as the vector of regression coefficients, this ridge regression has the following response vector $\bar{\mathbf{y}}$, design matrix $\bar{\mathbf{X}}$, and conditional ridge penalty matrix $\boldsymbol{\Omega}_{v|u}$:

$$\bar{\mathbf{y}} = \begin{pmatrix} \mathbf{x}_1 \\ \mathbf{x}_2 \\ ... \\ \mathbf{x}_m \end{pmatrix}, \qquad \bar{\mathbf{X}} = \begin{pmatrix} \mathbf{u} & \mathbf{0} & ... & \mathbf{0} \\ \mathbf{0} & \mathbf{u} & ... & \mathbf{0} \\ ... & ... & ... & ... \\ \mathbf{0} & \mathbf{0} & ... & \mathbf{u} \end{pmatrix},$$
$$\boldsymbol{\Omega}_{v|u} = \alpha \, \|\mathbf{u}\|^2 \, \boldsymbol{\Omega},$$

where $\mathbf{x}_j$ is the $j$-th column of $\mathbf{X}$, and where $\bar{\mathbf{y}}$ is of size $mn \times 1$ and $\bar{\mathbf{X}}$ is of size $mn \times m$. Both the design matrix $\bar{\mathbf{X}}$ and the ridge penalty depend on $\mathbf{u}$. It follows immediately that the penalized sum of squares (7) is equal to

$$\| \bar{\mathbf{y}} - \bar{\mathbf{X}}\mathbf{v} \|^2 \; + \; \mathbf{v}^T \, \boldsymbol{\Omega}_{v|u} \, \mathbf{v}, \tag{11}$$

which constitutes a penalized LS problem for $\mathbf{v}$. The associated penalized covariance matrix is

$$\bar{\mathbf{X}}^T \bar{\mathbf{X}} + \boldsymbol{\Omega}_{v|u} \; = \; (\mathbf{u}^T \mathbf{u}) \, (\mathbf{I} + \alpha \, \boldsymbol{\Omega}),$$

and thus the hat matrix of the ridge regression is

$$\mathbf{H} \; = \; \bar{\mathbf{X}} \, (\bar{\mathbf{X}}^T \bar{\mathbf{X}} + \boldsymbol{\Omega}_{v|u})^{-1} \, \bar{\mathbf{X}}^T \; = \; \frac{1}{\mathbf{u}^T \mathbf{u}} \, \bar{\mathbf{X}} \, \mathbf{S} \, \bar{\mathbf{X}}^T,$$

where $\mathbf{S} = \mathbf{S}(\alpha) = (\mathbf{I} + \alpha \, \boldsymbol{\Omega})^{-1}$.

Consider now the cross-validation that deletes one column of $\mathbf{X}$ at a time. This corresponds to deleting a block of size $n$ from $\bar{\mathbf{y}}$ at a time. Partition the hat matrix $\mathbf{H}$ into $m \times m$ equal-sized blocks where each block corresponds to a column of $\mathbf{X}$. Let $\hat{\mathbf{v}}^{(-j)} = (\hat{v}_1^{(-j)}, \ldots, \hat{v}_m^{(-j)})^T$ be the $\mathbf{v}$ that minimizes (11) when the $j$-th block of $\bar{\mathbf{y}}$ and the corresponding rows of $\bar{\mathbf{X}}$ are removed. We have the following lemma about the leave-out-one-column prediction errors.

**Lemma 1.** *The $j$-th leave-out-one-column cross-validated prediction error sum of squares is*

$$\|\mathbf{u}\hat{v}_j^{(-j)} - \mathbf{x}_j\|^2 = \mathbf{x}_j^T \mathbf{x}_j - \frac{(\mathbf{x}_j^T \mathbf{u})^2}{\|\mathbf{u}\|^2} + \frac{\left(\|\mathbf{u}\|\hat{v}_j - \frac{\mathbf{u}^T \mathbf{x}_j}{\|\mathbf{u}\|}\right)^2}{(1 - \gamma_j \|\mathbf{u}\|^2)^2}, \tag{12}$$

*where $\gamma_j = \mathbf{S}_{jj}/\|\mathbf{u}\|^2$, and $\mathbf{S}_{jj}$ is the $(j,j)$-th element of the matrix $\mathbf{S} = \mathbf{S}(\alpha)$.*

The proof of Lemma 1 is relegated to Appendix A. Note that $\hat{v}_j$ is the $j$-th element of $\hat{\mathbf{v}} = \mathbf{S} \, \mathbf{X}^T \mathbf{u}/\mathbf{u}^T \mathbf{u}$ and $\mathbf{x}_j^T \mathbf{u}$ is the $j$-th element of $\mathbf{X}^T \mathbf{u}$. Note also that $\gamma_j \|\mathbf{u}\|^2 = \mathbf{S}_{jj}$. Since we condition on $\mathbf{u}$, the first two terms on the right hand side of (12) are irrelevant. Averaging the last term in (12) over $j$, we obtain

$$\frac{1}{m} \sum_{j=1}^{m} \frac{\left\| \|\mathbf{u}\| \hat{v}_j - \frac{1}{\|\mathbf{u}\|} \mathbf{x}_j^T \mathbf{u} \right\|^2}{\left(1 - \{\mathbf{S}(\alpha)\}_{jj}\right)^2} = \frac{1}{m} \sum_{j=1}^{m} \frac{\left[\{(\mathbf{I} - \mathbf{S})(\mathbf{X}^T \mathbf{u})\}_{jj}\right]^2}{\|\mathbf{u}\|^2 \left(1 - \{\mathbf{S}(\alpha)\}_{jj}\right)^2}, \qquad (13)$$

which, ignoring the $\|\mathbf{u}\|^2$ factor in the denominator, is exactly the cross-validation criterion given in (9). Replacing $[\mathbf{S}(\alpha)]_{jj}$ in (13) by their average value, $(1/m)\mathrm{tr}\{\mathbf{S}(\alpha)\}$, we obtain

$$\frac{\frac{1}{m} \left\| \|\mathbf{u}\| \hat{\mathbf{v}} - \frac{1}{\|\mathbf{u}\|} \mathbf{X}^T \mathbf{u} \right\|^2}{\left(1 - \frac{1}{m}\mathrm{tr}\{\mathbf{S}(\alpha)\}\right)^2} = \frac{\frac{1}{m} \left\|(\mathbf{I} - \mathbf{S})(\mathbf{X}^T \mathbf{u})\right\|^2}{\|\mathbf{u}\|^2 \left(1 - \frac{1}{m}\mathrm{tr}\{\mathbf{S}(\alpha)\}\right)^2},$$

which, ignoring the $\|\mathbf{u}\|^2$ factor, is the generalized cross-validation criterion given in (10).

**Remark 2.** Both Rice and Silverman (1991) and Silverman (1996) suggest cross-validation based on row deletion in $\mathbf{X}$. Row deletion lacks simple computational shortcuts such as (9); it hence involves actual computation of a large number of leave-out-one-row estimates. As a result, CV based on row deletion is computationally expensive as confirmed by our numerical studies in Section 6.

**Remark 3.** Since the principal components are extracted sequentially, the squared errors in the numerators of the CV and GCV criteria vary with the principal components. The CV score (9) and the GCV score (10) are defined for extracting the first principal component. When each subsequent principal component is extracted, the $\mathbf{X}$ matrix in the definition of the CV and GCV scores needs to be replaced by the residual matrix after the effects of the previous principal components are removed.

## 5. Extracting the whole curve of principal component weight function

So far we have focused on the discretized problem, although the functional nature has been taken into account by regularization with a second-order roughness penalty. This section introduces an optimization criterion, the minimizer of which gives the entire principal component weight function. It also turns out that the extracted function is a natural cubic spline that interpolates the weighted vector obtained by minimizing (7). Our development relies on some standard results from spline smoothing that can be found in Green and Silverman (1994).

Replacing the discretized vector $\mathbf{v}$ with the complete function $\gamma(\cdot)$, we propose to find the estimate $\hat{\gamma}(\cdot)$ of the first principal component weight function by minimizing, with respect to $u_i$ and $\gamma(\cdot)$, the penalized sum of squares,

$$\sum_{i=1}^{n} \sum_{j=1}^{m} \{x_{ij} - u_i \gamma(t_j)\}^2 + \alpha \left(\sum_{i=1}^{n} u_i^2\right) \int \{\gamma''(t)\}^2 \, dt, \qquad (14)$$

where $\alpha > 0$ is a smoothing parameter. The estimated principal component function $\hat{\gamma}(\cdot)$ is the optimizer of (14) over the class of all functions that satisfy $\int \{\gamma''(t)\}^2\, dt < \infty$. Similar to (7), the two terms in (14) have the same units, and the optimization problem is invariant under scale transformation of the measurements. However, unlike (7), the whole function is recovered by optimizing (14).

We now characterize the solution to the optimization problem (14) and show that it is closely related to the solution of the discretized problem (7) with an appropriately defined penalty matrix. To this end, we define two banded matrices, $Q$ and $R$, as follows. Let $h_j = t_{j+1} - t_j$ for $j = 1, \ldots, m-1$. Let $Q$ be the $m \times (m-2)$ matrix with entries $q_{jk}$, for $j = 1, \ldots, m$ and $k = 2, \ldots, m-1$, given by

$$q_{k-1,k} = h_{k-1}^{-1},\ q_{kk} = -h_{k-1}^{-1} - h_k^{-1},\ q_{k+1,k} = h_k^{-1}$$

for $k = 2, \ldots, m-1$, and $q_{jk} = 0$ for $|j - k| > 2$. To simplify the presentation, the columns of $Q$ are numbered in a non-standard way, starting at $k = 2$, so that the top left element of $Q$ is $q_{12}$. The symmetric matrix $R$ is $(m-2) \times (m-2)$ with elements $r_{jk}$, for $j$ and $k$ running from 2 to $(m-1)$, given by

$$r_{jj} = \frac{1}{3}(h_{j-1} + h_j) \text{ for } j = 2, \ldots, m-1,$$

$$r_{j,j+1} = r_{j+1,j} = \frac{1}{6}h_j \text{ for } j = 2, \ldots, m-2,$$

and $r_{jk} = 0$ for $|j - k| > 2$. The matrix $R$ is strictly diagonal dominant and thus is strictly positive-definite.

**Theorem 1.** *The $\hat{\gamma}(\cdot)$ optimizing (14) is a natural cubic spline with knots at $t_j$. Let $\hat{v}_j = \hat{\gamma}(t_j)$ and $\hat{\mathbf{v}} = (\hat{v}_1, \ldots, \hat{v}_m)^T$. Then $\hat{\mathbf{v}}$ is the optimizer of the discretized problem (7) with the penalty matrix $\boldsymbol{\Omega} = QR^{-1}Q^T$.*

The proof of Theorem 1 can be found in Appendix B. According to Theorem 1, to obtain the entire curve of the principal component weight function, one needs to first solve (7) with a penalty matrix $\boldsymbol{\Omega} = QR^{-1}Q^T$ to obtain $\hat{v}_1, \ldots, \hat{v}_m$, and then find the natural cubic spline that interpolates $(t_j, \hat{v}_j)$. Computation of the interpolating natural cubic spline $\gamma(\cdot)$ at any evaluation point $t$ using its values at the knots $t_j$ is a standard operation. Specifically, let $v_j = \gamma(t_j)$ and $s_j = \gamma''(t_j)$. By the definition of natural cubic spline, the second derivative of $\gamma$ at $t_1$ and $t_m$ is zero, so that $s_1 = s_m = 0$. The interpolating natural cubic spline is completely determined by its values and second derivatives at each of the knots $t_j$ according to the following formula:

$$\gamma(t) = \frac{(t - t_j)v_{j+1} + (t_{j+1} - t)v_j}{h_j}$$
$$- \frac{1}{6}(t - t_j)(t_{j+1} - t)\left\{\left(1 + \frac{t - t_j}{h_j}\right)s_{j+1} + \left(1 + \frac{t_{j+1} - t}{h_j}\right)s_j\right\}$$

for $t_j \leq t \leq t_{j+1}$, $j = 1, \ldots, m-1$; the spline outside $[t_1, t_m]$ is obtained by linear extrapolation. The vector $\mathbf{s} = (s_2, \ldots, s_{m-1})^T$ of the second derivatives used above can be obtained by $\mathbf{s} = R^{-1}Q^T\mathbf{v}$. See Chapter 2 of Green and Silverman (1994) for more details of efficient computation of interpolating splines.

# 6. Numerical results

## 6.1. Call center arrival data example

We applied the proposed method to the call center arrival data analyzed in Shen and Huang (2008). The data recorded the number of calls that got connected to a call center during every quarter hour between 7:00AM and midnight for every weekday between January 1 and October 26 in the year 2003. In total, there are 42 whole weeks during the period and each day consists of 68 quarter hours. Let $N_{ij}$ denote the call volume during the $j$-th time interval on day $i$. We used the transformed data $X_{ij} = \sqrt{N_{ij} + 1/4}$ which together form a $210 \times 68$ matrix. The square-root transformation is used to stabilize variance and make the distribution close to normal. The same transformation has been used previously by Brown et al. (2005) and Shen and Huang (2008).

The mean curve from the data (or the column mean vector of $\mathbf{X}$) is quite smooth and summarizes the average intraday arrival pattern (Figure 1). It is bimodal with a main peak around 11:00AM followed by a second lower peak around 2:00PM. We subtracted the mean curve from the data and then ap-plied the penalized sum of squares criterion (14) to sequentially extract the first three principal components. The leave-out-one-column CV is used to select the smoothing parameter for each principal component respectively. A set of grid points of $\alpha = 0$ and $\alpha = 1.5^i$ for $i = -5, \ldots, 25$ are examined as candidate values. The chosen values are 0.44, 25.63 and 38.44 respectively. The estimated smooth principal component weight functions are plotted in Figure 1. The application of the leave-out-one-column GCV leads to qualitatively very similar results.

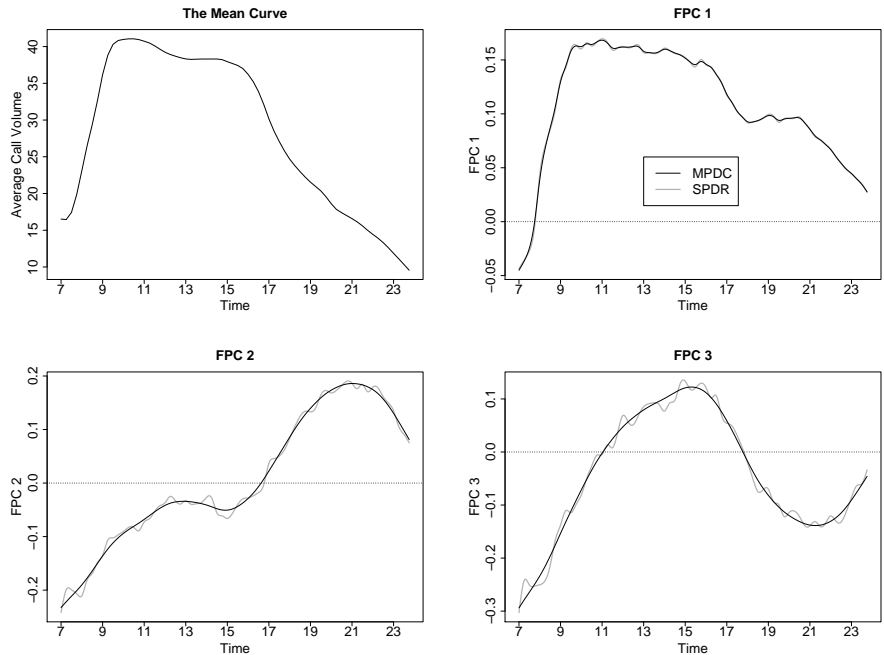

FIG 1. *The call center data. The mean curve and the estimated first three functional principal components. MPDC refers to our method of selecting multiple smoothing parameters using delete-column CV. SPDR refers to Silverman's method of selecting the single smoothing parameters using delete-row CV.*

For comparison, we also followed Silverman (1996)and applied CV with row deletion to select a single smoothing parameter. This alternative method selects $\alpha = 0.13$. The estimated principal component functions are plotted in Figure 1. We observe that Silverman's method undersmooths the principal component functions, while our method of using multiple smoothing parameters performs suitable smoothing. The computing times also demonstrate that the CV with row deletion is computationally much more expensive than our delete-column CV because of its lack of a computational shortcut. To generate the results in Figure 1, it took 55 and 2 seconds for the two approaches, respectively, using our R program running on a Debian Linux desktop with Intel® Pentium® 4 CPU of a clock speed of 2.8 Gigahertz. We used the computational tricks presented in Appendix C when implementing both approaches.

### 6.2. A simulated data example

To further understand the difference between the two methods used in Section 6.1, we performed the following simulation study. The data generating model is

$$X_{ij} = u_{i1}v_1(t_j) + u_{i2}v_2(t_j) + \epsilon_{ij}, \quad i = 1, \ldots, n; \; j = 1, \ldots, m, \quad (15)$$

where $u_{i1} \overset{i.i.d.}{\sim} N(0, \sigma_1^2)$, $u_{i2} \overset{i.i.d.}{\sim} N(0, \sigma_2^2)$, and $\epsilon_{ij} \overset{i.i.d.}{\sim} N(0, \sigma^2)$. The parameters are chosen as $n = m = 101$, $\sigma_1 = 20$, $\sigma_2 = 10$, $\sigma = 4$, and the 101 grid points $t_j$ are equally spaced in $[-1, 1]$. The two underlying functional principal component functions are

$$v_1(t) = \frac{1}{s_1}\{t + \sin(\pi t)\} \quad \text{and} \quad v_2(t) = \frac{1}{s_2}\cos(3\pi t),$$

where $s_1$ and $s_2$ are the normalizing constants that ensure $v_1$ and $v_2$ to have unit norm. We generated one hundred simulated data sets, and estimated the first two smooth principal component functions for each data set. We calculated the mean squared errors (MSE) over the 101 grid points for each estimated principal component function.

Panels (a) and (b) of Figure 2 are the scatterplots of MSEs for the two methods, with each point representing a simulated data set. We observe that using single smoothing parameter yields larger MSEs for most simulated data sets, especially for the first functional principal component (FPC). The same message is also evident in Panels (c) and (d) that provide the histograms of the ratios of MSEs of the two methods. Table 1 reports some summary statistics of the MSE ratio and shows that, by using a single instead of multiple smoothing parameters, the mean of the MSE ratio increases by 64% for the first FPC, and by 8% for the second FPC. We confirm the observed difference as significant by the sign test that gives the *p*-values that are essentially 0 for both FPCs. We also observe that the difference of the two methods for estimating the second FPC is not as big as that for the first one. This is because the two methods select similar smoothing parameters for the second FPC. As far as computing times go, delete-row CV on average needs 116 seconds for one simulation while delete-column CV only takes 2 seconds, using the same computer reported in Section 6.1.

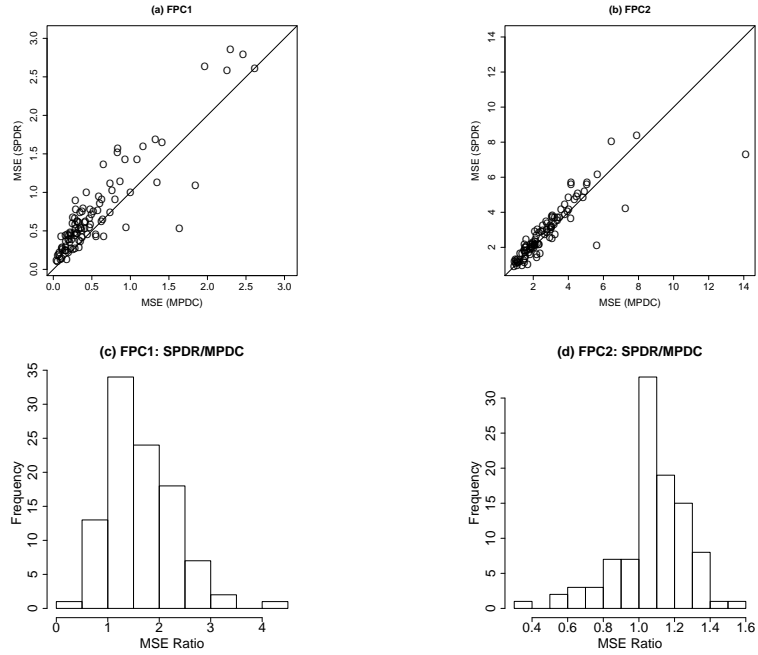

Fig 2. *Comparison of two methods for estimating functional principal components by simulation study. MPDC refers to our method of selecting multiple smoothing parameters using delete-column CV. SPDR refers to Silverman's method of selecting the single smoothing parameters using delete-row CV. The reported MSEs equal to the actual values multiplied by $10^4$.*

Table 1
*Summary statistics of the MSE ratios for the two methods in Figure 2. Q1 and Q3 refer to the lower and upper quartiles.*

| | FPC 1 | | | | FPC 2 | | | |
|---|---|---|---|---|---|---|---|---|
| | $Q1$ | Median | Mean | $Q3$ | $Q1$ | Median | Mean | $Q3$ |
| SPDR/MPDC | 1.17 | 1.51 | 1.64 | 2.03 | 1.01 | 1.08 | 1.07 | 1.20 |

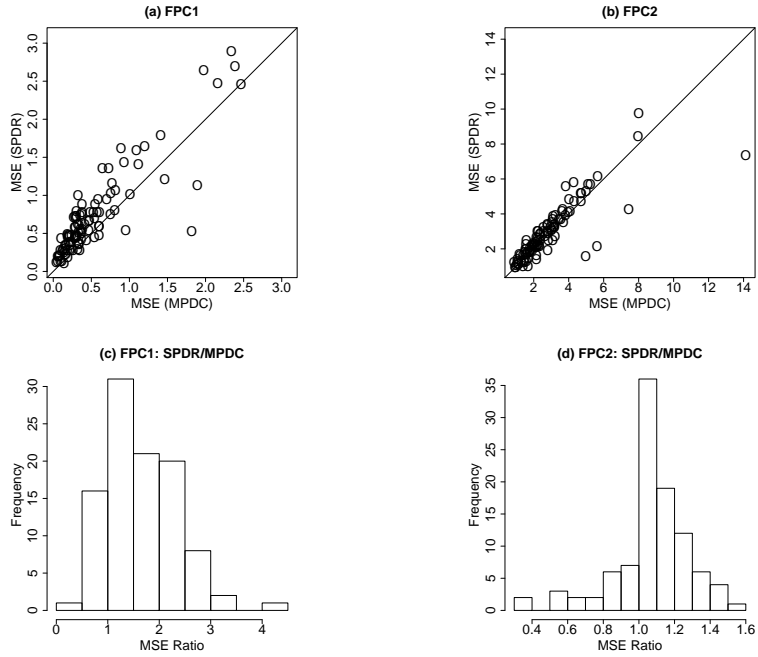

Fɪɢ 3. *Same as Figure 2, except that a mean function has to be estimated.*

TABLE 2
*Same as Table 1, except that a mean function has to be estimated.*

|  | FPC 1 | | | | FPC 2 | | | |
|---|---|---|---|---|---|---|---|---|
|  | *Q*1 | Median | Mean | *Q*3 | *Q*1 | Median | Mean | *Q*3 |
| SPDR/MPDC | 1.22 | 1.53 | 1.66 | 2.08 | 1.01 | 1.08 | 1.07 | 1.19 |

To see the effect of estimating the mean function on extracting the principal component functions, a mean function was added to the data generating model (15). The mean was then removed by column centering of the data matrix prior to applying the FPCA algorithms. Figure 3 and Table 2 present the results after removing the mean, and when compared with Figure 2 and Table 1, suggest that the effect of estimating the mean function is not very significant.

## Appendix A: Proof of Lemma 1

Note that $\hat{\mathbf{v}}^{(-j)}$ also solves the ridge regression (11) when the $j$-th block of $\bar{\mathbf{y}}$ is replaced by $\mathbf{u}\hat{v}_j^{(-j)}$. The $j$-th block of the fitted equation $\hat{\bar{\mathbf{y}}} = \mathbf{H}\bar{\mathbf{y}}$ of this latter ridge regression reads as

$$\mathbf{u}\hat{v}_j^{(-j)} = \sum_{k \neq j} \mathbf{H}_{jk}\mathbf{x}_k + \mathbf{H}_{jj}\{\mathbf{u}\hat{v}_j^{(-j)}\}.$$

Subtracting $\mathbf{x}_j$ on both sides of the above identity and observing that $\sum_k \mathbf{H}_{jk}\mathbf{x}_k = \mathbf{u}\hat{\mathbf{v}}_j$, we obtain

$$\mathbf{u}\hat{v}_j^{(-j)} - \mathbf{x}_j = \sum_k \mathbf{H}_{jk}\mathbf{x}_k - \mathbf{x}_j + \mathbf{H}_{jj}\{\mathbf{u}\hat{v}_j^{(-j)} - \mathbf{x}_j\}$$

$$= \mathbf{u}\hat{v}_j - \mathbf{x}_j + \mathbf{H}_{jj}\{\mathbf{u}\hat{v}_j^{(-j)} - \mathbf{x}_j\}.$$

Therefore, the cross-validated residual for deleting the $j$-th column of $\mathbf{X}$ is

$$\mathbf{u}\hat{v}_j^{(-j)} - \mathbf{x}_j = (\mathbf{I} - \mathbf{H}_{jj})^{-1}(\mathbf{u}\hat{v}_j - \mathbf{x}_j)$$

where

$$\mathbf{H}_{jj} = \frac{\mathbf{S}_{jj}}{\mathbf{u}^T\mathbf{u}}\mathbf{u}\mathbf{u}^T = \gamma_j\mathbf{u}\mathbf{u}^T.$$

Denote $\mathbf{w} = \mathbf{u}\hat{v}_j - \mathbf{x}_j$. Its squared norm is

$$\|\mathbf{w}\|^2 = \mathbf{x}_j^T\mathbf{x}_j - 2\mathbf{x}_j^T\mathbf{u}\hat{v}_j + \|\mathbf{u}\|^2\hat{v}_j^2$$

$$= \mathbf{x}_j^T\mathbf{x}_j - \frac{(\mathbf{x}_j^T\mathbf{u})^2}{\|\mathbf{u}\|^2} + \left(\|\mathbf{u}\|\hat{v}_j - \frac{\mathbf{u}^T\mathbf{x}_j}{\|\mathbf{u}\|}\right)^2. \tag{A.1}$$

Since $\mathbf{u}^T\mathbf{w} = \|\mathbf{u}\|^2(\hat{v}_j - \mathbf{u}^T\mathbf{x}_j/\|\mathbf{u}\|^2)$, we have that

$$\frac{(\mathbf{u}^T\mathbf{w})^2}{\|\mathbf{u}\|^2} = \left(\|\mathbf{u}\|\hat{v}_j - \frac{\mathbf{u}^T\mathbf{x}_j}{\|\mathbf{u}\|}\right)^2. \tag{A.2}$$

Using the identity

$$(\mathbf{I} - \gamma_j\mathbf{u}\mathbf{u}^T)^{-1} = \mathbf{I} + \frac{\gamma_j}{1 - \gamma_j\|\mathbf{u}\|^2}\mathbf{u}\mathbf{u}^T,$$

we can write the cross-validated residual $\mathbf{u}\hat{v}_j^{(-j)} - \mathbf{x}_j$ as

$$(\mathbf{I} - \mathbf{H}_{jj})^{-1}(\mathbf{u}\hat{v}_j - \mathbf{x}_j) = \left(\mathbf{I} + \frac{\gamma_j}{1 - \gamma_j\|\mathbf{u}\|^2}\mathbf{u}\mathbf{u}^T\right)\mathbf{w}$$

$$= \mathbf{w} + \frac{\gamma_j}{1 - \gamma_j\|\mathbf{u}\|^2}(\mathbf{u}^T\mathbf{w})\mathbf{u}.$$

Thus the squared norm of $\mathbf{u}\hat{v}_j^{(-j)} - \mathbf{x}_j$ is

$$\|\mathbf{w}\|^2 + \frac{2\gamma_j}{1 - \gamma_j\|\mathbf{u}\|^2}(\mathbf{u}^T\mathbf{w})^2 + \frac{\gamma_j^2}{(1 - \gamma_j\|\mathbf{u}\|^2)^2}(\mathbf{u}^T\mathbf{w})^2\|\mathbf{u}\|^2$$

$$= \|\mathbf{w}\|^2 + \frac{(\mathbf{u}^T\mathbf{w})^2}{\|\mathbf{u}\|^2}\left\{\frac{2\gamma_j\|\mathbf{u}\|^2}{(1 - \gamma_j\|\mathbf{u}\|^2)} + \frac{\gamma_j^2\|\mathbf{u}\|^4}{(1 - \gamma_j\|\mathbf{u}\|^2)^2}\right\}$$

$$= \|\mathbf{w}\|^2 + \frac{(\mathbf{u}^T\mathbf{w})^2}{\|\mathbf{u}\|^2}\left\{\frac{1}{(1 - \gamma_j\|\mathbf{u}\|^2)^2} - 1\right\}.$$

Combining this result with (A.1) and (A.2) we obtain (12). □

## Appendix B: Proof of Theorem 1

Since the natural cubic spline interpolant is the unique minimizer of $\int \gamma''^2$ over all functions that interpolate the data $(t_j, \hat{v}_j)$ (Theorem 2.3 of Green and Silverman, 1994), the minimizing function $\hat{\gamma}(\cdot)$ of (14) is necessarily a natural cubic spline with knots at the points $t_j$. Therefore in (14) we can restrict attention to natural cubic splines with knots at $t_j$. A natural cubic spline $\gamma(t)$ that interpolates $(t_j, v_j)$ is uniquely defined. By Theorem 2.1 of Green and Silverman (1994), $\int \{\gamma''(t)\}^2 \, dt = \mathbf{v}^T \boldsymbol{\Omega} \mathbf{v}$, where $\mathbf{v} = (v_1, \dots, v_m)^T$. Let $\mathbf{X} = (x_{ij})$ and $\mathbf{u} = (u_1, \dots, u_n)^T$. It follows immediately that the penalized sum of squares (14) can be written as

$$\|\mathbf{X} - \mathbf{u}\mathbf{v}^T\|_F^2 + \alpha \, \mathbf{u}^T \mathbf{u} \, \mathbf{v}^T \boldsymbol{\Omega} \mathbf{v},$$

which is exactly the discretized criterion (7) given in Section 2. $\qquad \square$

## Appendix C: Implementation details

We provide in this appendix some implementation details of our method. We first describe a method of computing smoothed PC using any existing SVD software, then discuss efficient computation of the CV/GCV criteria for multiple candidate values of $\alpha$.

Note that the penalized reconstruction error criterion (7) can be expanded as follows:

$$\|\mathbf{X} - \mathbf{u}\mathbf{v}^T\|_F^2 + \alpha \, \mathbf{u}^T \mathbf{u} \, \mathbf{v}^T \boldsymbol{\Omega} \mathbf{v} = \|\mathbf{X}\|_F^2 - 2\mathbf{u}^T \mathbf{X} \mathbf{v} + \mathbf{u}^T \mathbf{u} \, \mathbf{v}^T (\mathbf{I} + \alpha \, \boldsymbol{\Omega}) \mathbf{v}.$$

Denote $\mathbf{S}(\alpha) = (\mathbf{I} + \alpha \, \boldsymbol{\Omega})^{-1}$, $\tilde{\mathbf{v}} = \mathbf{S}^{-1/2}(\alpha) \, \mathbf{v}$, and $\tilde{\mathbf{X}} = \mathbf{X} \, \mathbf{S}^{1/2}(\alpha)$. Then, the above expression is equivalent to

$$\|\mathbf{X}\|_F^2 - \|\tilde{\mathbf{X}}\|_F^2 + \|\tilde{\mathbf{X}}\|_F^2 - 2\mathbf{u}^T \, \tilde{\mathbf{X}} \, \tilde{\mathbf{v}} + \mathbf{u}^T \mathbf{u} \, \tilde{\mathbf{v}}^T \tilde{\mathbf{v}},$$

which can be simplified as

$$\|\mathbf{X}\|_F^2 - \|\tilde{\mathbf{X}}\|_F^2 + \|\tilde{\mathbf{X}} - \mathbf{u}\tilde{\mathbf{v}}^T\|_F^2. \tag{A.3}$$

For a given smoothing parameter $\alpha$, the minimizing $\mathbf{u}$ and $\tilde{\mathbf{v}}$ of (A.3) can be easily obtained as the first pair of singular vectors of $\tilde{\mathbf{X}} = \mathbf{X}\mathbf{S}^{1/2}(\alpha)$, as discussed in Section 2.1. Since $\mathbf{S}^{1/2}(\alpha)$ can be interpreted as a half-smoothing operator, the transformed matrix $\tilde{\mathbf{X}}$ is obtained by half-smoothing the rows of the original data matrix $\mathbf{X}$. After $\tilde{\mathbf{v}}$ is obtained as the first right singular vector of $\tilde{\mathbf{X}}$, we half-smooth it to obtain the smoothed PC function $\mathbf{v} = \mathbf{S}^{1/2}(\alpha)\tilde{\mathbf{v}}$. Thus by using existing SVD software, we can avoid directly programming the iterative power algorithm. This is convenient when a high level programming language is used for coding.

The eigen decomposition of the symmetric and positive definte penalty matrix $\boldsymbol{\Omega}$ can be written as

$$\boldsymbol{\Omega} = \boldsymbol{\Gamma} \boldsymbol{\Lambda} \boldsymbol{\Gamma}^T,$$

where $\mathbf{\Gamma}$ is an orthogonal matrix containing the eigenvectors, and $\mathbf{\Lambda}$ is a diagonal matrix of the eigenvalues. For any value of $\alpha$, the eigen decompositions of the smoothing and half-smoothing matrices are

$$\mathbf{S}(\alpha) = \mathbf{\Gamma}\,(\mathbf{I} + \alpha\,\mathbf{\Lambda})^{-1}\mathbf{\Gamma}^T \quad \text{and} \quad \mathbf{S}^{1/2}(\alpha) = \mathbf{\Gamma}\,(\mathbf{I} + \alpha\,\mathbf{\Lambda})^{-1/2}\mathbf{\Gamma}^T.$$

Using the eigen decomposition of $\mathbf{S}^{1/2}(\alpha)$, we see that minimization of $\|\tilde{\mathbf{X}} - \mathbf{u}\tilde{\mathbf{v}}^T\|_F^2$ is equivalent to minimization of $\|\mathbf{X}\mathbf{\Gamma}(\mathbf{I} + \alpha\mathbf{\Lambda})^{-1/2} - \mathbf{u}\bar{\mathbf{v}}^T\|_F^2$ where $\bar{\mathbf{v}} = \mathbf{\Gamma}^T\tilde{\mathbf{v}}$. Thus, after obtaining the first right singular vector $\bar{\mathbf{v}}$ of $\mathbf{X}\mathbf{\Gamma}(\mathbf{I} + \alpha\mathbf{\Lambda})^{-1/2}$, we obtain $\mathbf{v}$ using $\mathbf{v} = \mathbf{S}^{1/2}(\alpha)\tilde{\mathbf{v}} = \mathbf{\Gamma}(\mathbf{I} + \alpha\mathbf{\Lambda})^{-1/2}\bar{\mathbf{v}}$. Note that $\mathbf{X}\mathbf{\Gamma}$ only needs to be computed once when $\mathbf{v}$ needs to be computed for a set of values of $\alpha$.

Now we discuss efficient evaluation of the GCV criterion (10). Since

$$\mathbf{I} - \mathbf{S}(\alpha) = \mathbf{\Gamma}\{\mathbf{I} - (\mathbf{I} + \alpha\mathbf{\Lambda})^{-1}\}\mathbf{\Gamma}^T,$$

the trace that appears in the denominator of (10) equals to

$$\text{tr}\{\mathbf{S}(\alpha)\} = \sum_k \frac{1}{1 + \alpha\lambda_k},$$

and the numerator in (10) is

$$\|\{\mathbf{I} - \mathbf{S}(\alpha)\}(\mathbf{X}^T\mathbf{u})\| = \|\{\mathbf{I} - (\mathbf{I} + \alpha\mathbf{\Lambda})^{-1}\}(\mathbf{X}\mathbf{\Gamma})^T\mathbf{u}\|.$$

Denote $\mathbf{w} = (\mathbf{X}\mathbf{\Gamma})^T\mathbf{u}$. The numerator of the GCV criterion equals to the Euclidean norm of the shrunken $\mathbf{w}$ with the $k$-th component shrunken by a factor of $\alpha\lambda_k/(1 + \alpha\lambda_k)$. When computation of GCV is needed for a set of candidate values of $\alpha$, considerable computing saving is achieved since $\mathbf{X}\mathbf{\Gamma}$ and the eigen decomposition of $\mathbf{\Omega}$ only need to be calculated once.